\documentclass[a4paper,12pt]{amsart}
\usepackage{amsfonts}
\usepackage{mathrsfs}
\usepackage{amsmath,amssymb,latexsym,amsfonts,amscd}

\title{Characterization of the 4-canonical birationality of algebraic
threefolds}
\author{Meng Chen and De-Qi Zhang}
\address{\rm Institute of Mathematics, School of Mathematical Sciences,
Fudan University, Shanghai, 200433, China}
\email{mchen@fudan.edu.cn}

\address{\rm Department of Mathematics, National University of Singapore,
2 Science Drive 2, Singapore 117543, Republic of Singapore}
\email{matzdq@nus.edu.sg}

\thanks{The first author was supported by both the Program for New Century
Excellent Talents in University (\#NCET-05-0358) and the
National Outstanding Young Scientist Foundation (\#10625103).
The second author was partially supported by an Academic
Research Fund of NUS }

\newcommand{\OO}{{\mathcal{O}}}

\newcommand{\bQ}{{\mathbb Q}}
\newcommand{\bP}{{\mathbb P}}
\newcommand{\roundup}[1]{\ulcorner{#1}\urcorner}

\newtheorem{thm}{Theorem}[section]
\newtheorem{lem}[thm]{Lemma}

\newtheorem{prop}[thm]{Proposition}
\newtheorem{claim}[thm]{Claim}
\theoremstyle{definition}
\newtheorem{defn}[thm]{Definition}
\newtheorem{setup}[thm]{}

\newtheorem{exmp}[thm]{Example}

\theoremstyle{remark}

\begin{document}


\begin{abstract}
In this article we present a 3-dimensional analogue of a well-known
theorem of E. Bombieri (in 1973) which characterizes the
bi-canonical birationality of surfaces of general type. Let $X$ be a
projective minimal 3-fold of general type with
$\mathbb{Q}$-factorial terminal singularities and the geometric
genus $p_g(X)\ge 5$. We show that the 4-canonical map $\varphi_4$ is
{\it not} birational onto its image if and only if $X$ is
birationally fibred by a family $\mathscr{C}$ of irreducible curves
of geometric genus 2 with $K_X\cdot C_0=1$ where $C_0$ is a general
irreducible member in $\mathscr{C}$.
\end{abstract}
\maketitle
\pagestyle{myheadings} \markboth{\hfill Meng Chen and De-Qi
Zhang\hfill}{\hfill Characterization of the 4-canonical
birationality\hfill}
\section{\bf Introduction}
We work over the field $\mathbb{C}$ of complex numbers.

In this article we study pluri-canonical systems of projective minimal
threefolds $X$ of general type.

Our guidance are results of E. Bombieri \cite{Bom} on minimal
surfaces $S$. For $m\geq 5$, the $m$-canonical map $\Phi_{|mK_S|}$
is always birational. If, for $1<m<5$, $\Phi_{|mK_S|}$ is not
birational and if the numerical invariants of $S$ are sufficiantly
large, then $S$ carries a pencil of curves of small genus. For
example Bombieri shows the following (see \cite{Bom}, or \cite{BPV},
Theorem 5.4(iv)):

\begin{thm} Let $S$ be a smooth minimal projective surface of
general type with $K_S^2\ge 10$ and $p_g(S)\ge 6$. Then the
bi-canonical map is not birational onto its image if and only if $S$
has a family of curves of genus $2$.
\end{thm}

For minimal threefolds $X$ with $\mathbb{Q}$-factorial terminal
singularities it is known that there exists some universal
constant $r_3$ such that the pluri-canonical map
$\varphi_{r_3}:=\Phi_{|r_3K_X|}$ is birational (see Tsuji
\cite{Tsuji}, Hacon-M$^{\rm c}$Kernan \cite{H-M} and Takayama
\cite{Tak}). Tsuji \cite{Tsuji} has ever proved $r_3\leq
18(2^9\cdot 3^7)!$. If one requires in addition that either the
invariants of $X$ are big (e.g. for $p_g(X)\ge 4$ see
\cite{IJM}; for $K_X^3\gg 0$ see \cite{Tod}) or $X$ is
Gorenstein (see \cite{Crelle}), one may take $r_3=5$. We remark
that $r_3$ can not be 4 because the 4-canonical map of the
product of a curve and a surface of type $(K^2, p_g) = (1,2)$ is not
birational. So it is natural to ask how $\Phi_{|4K_X|}$ behaves,
provided again that the numerical invariants of $X$ are large.

\begin{setup}{\bf Known works.} There
are several works that deal with the behaviour of $\varphi_4$.
In 2000, S. Lee (\cite{Lee2}) proved the base point freeness of
$|4K|$ for minimal Gorenstein threefolds of general type. In
2002, the first author (\cite{Tianjin}) gave a sufficient
condition for the birationality of $\varphi_4$ again for minimal
Gorenstein threefolds. Then, in 2005, J. Dong (\cite{D})
improved upon the method of the first author and gave a
sufficient condition for the birationality of $\varphi_4$ for
arbitrary threefolds of general type.
\end{setup}

Our main result is the following:

\begin{thm}\label{main} Let $X$ be a minimal projective threefold of
general type with $\mathbb{Q}$-factorial terminal singularities and
the geometric genus $p_g(X)\ge 5$. Then:
\begin{itemize}
\item
[(i)] The $4$-canonical map $\varphi_4:= \Phi_{|4K_X|} $ is
generically finite of degree $\le 2$.
\item
[(ii)] $\varphi_4$ is {\bf not} birational if and only if $X$ is
birationally fibred by a family $\mathscr{C}$ of irreducible curves
of geometric genus $2$ with $K_X\cdot C_0=1$ for a general member
$C_0$ in $\mathscr{C}$.
\item
[(iii)] In \rm{(ii)} the family $\mathscr{C}$ is birationally, uniquely
determined by the given threefold $X$.
\end{itemize}
\end{thm}

The precise definition of ``birationally fibred by a family
$\mathscr{C}$ of irreducible curves'' will be given in Definition
\ref{coverf}. We remark that the condition $p_g(X) \ge 5$ in Theorem
\ref{main} is optimal (see Example \ref{4}).

If a minimal 3-fold $X$ is birationally fibred by surfaces of type
$(K^2, p_g)$ $= (1,2)$, then $\varphi_4$ is not birational.  Here we mean that there
exists a birational morphism $\nu:Z\longrightarrow X$ and a
fibration $\tilde{f}: Z\longrightarrow \tilde{B}$ with $Z$ a smooth threefold
and $\tilde{B}$ a smooth curve, such that a general fiber of
$\tilde{f}$ is a surface of type (1,2). Taking further birational
modification of $Z$ if necessary, we may assume that the relative
canonical map $\tilde{\Psi}: Z\longrightarrow
\mathbb{P}(\tilde{f}_*\omega_{Z/\tilde{B}}^\vee)$ is a morphism over
$\tilde{B}$.

Our next result shows how to find the unique family $\mathscr{C}$ of
curves mentioned in Theorem \ref{main} from the family of $(1, 2)$
surfaces on $X$.

\begin{thm}\label{(1,2)}
Let $X$ be a minimal projective threefold of general type with
$\mathbb{Q}$-factorial terminal singularities and the geometric
genus $p_g(X)\ge 5$. Suppose that $X$ is birationally fibred by
surfaces of type $(K^2, p_g) = (1, 2)$. Then the birationally unique family
$\mathscr{C}$ in {\rm Theorem \ref{main} (ii)} is on the fibres of
$\tilde{f}$ and induced by $\tilde{\Psi}$.
\end{thm}

Theorem \ref{main} has the following application to the effect that
minimal threefolds $X$ of general type with small slope $(K^3/p_g)$
will have non-birational $4$-canonical maps $\varphi_4$.

\begin{thm}\label{app} Let $X$ be a minimal projective threefold of
general type with
$\mathbb{Q}$-factorial terminal singularities. Then $\varphi_4$ is
not birational whenever either of the following two conditions is
satisfied:
\begin{itemize}
\item
[(i)] $K_X^3<\frac{4}{3}(p_g(X)-2)$ and $p_g(X)\not\in [3, 11] $.
\item
[(ii)] $X$ is Gorenstein and $K_X^3<2p_g(X)-6$.
\end{itemize}
\end{thm}

\par
There are infinitely many non-trivial examples satisfying the
conditions in Theorem \ref{app} (see the last section).

\section{\bf Notation and the set up}
We use the standard notation and terminology in the textbook of
Hartshorne. Throughout the paper $D_1 \sim D_2 $ (resp. $=_{\bQ}$,
or $D_1 \equiv D_2$) means that divisors $D_1$ and $D_2$ are
linearly equivalent (resp. $rD_1$ and $rD_2$ are linearly equivalent
for some positive integer $r$, or $D_1$ and $D_2$ are numerically
equivalent). By {\it a minimal variety}, we always mean one with a
nef canonical divisor $K$ and with terminal singularities. Given a
smooth projective threefold $V$ of general type, the 3-dimensional
MMP (see \cite{K-M, KMM}) for example) says that $V$ has a minimal
model $X$ with ${\mathbb Q}$FT singularities. Since the
birationality of $\Phi_m$ is equivalent to that of $\varphi_m$, we
may begin our study from a minimal threefold $X$.

\begin{defn}\label{coverf} Let $Y$ be a $\bQ$-Gorenstein
(i.e., $r K_Y$ being Cartier for some positive integer $r$) normal
projective variety. We say that $Y$ is {\it birationally fibred by a
family of curves} if there exist both a birational morphism $\pi:
Y'\longrightarrow Y$ and a fibration $f_{Y'}:Y'\longrightarrow W$
whose general fiber is a smooth projective curve, where $Y'$ and $W$
are projective normal varieties. Denote by $\mathscr{C}$ the set
$\{\pi(Y'_w) \ | \ w\in W,\ Y'_w = f_{Y'}^{-1}(w)\}$. An element
$C_0\in \mathscr{C}$ is called {\em a generic member} if $C_0$ is
the $\pi$-image of a general fiber $Y_w' = f_{Y'}^{-1}(w)$ of
$f_{Y'}$. So the intersection number $(K_Y\cdot C_0) = (\pi^*K_Y
\cdot Y_w')$ is uniquely determined by the curve family
$\mathscr{C}$ and is independent of the choice of the birational
modification $\pi$.
\end{defn}

\begin{setup}\label{1-map}{\bf The canonical map.}
Suppose $p_g(X)\ge 2$. We may study the canonical map $\varphi_{1}$
which is only a rational map. First we fix an effective Weil divisor
$K_1\sim K_X$. Take successive blow-ups $\pi: X'\rightarrow X$
(along nonsingular centers), which exists by Hironaka's big theorem,
such that:
\par \vskip 1pc
(i) $X'$ is smooth;

(ii) the movable part of $|K_{X'}|$ is base point free;

(iii) the support of $\pi^*(K_1)$ is of simple normal crossings.
\par \vskip 1pc
Denote by $g$ the composition $\varphi_{1}\circ\pi$. So $g:
X'\longrightarrow W'\subseteq{\mathbb P}^{p_g(X)-1}$ is a morphism.
Let $X'\overset{f}\longrightarrow B\overset{s}\longrightarrow W'$ be
the Stein factorization of $g$. We have the following commutative
diagram:\medskip

\begin{picture}(50,80) \put(100,0){$X$} \put(100,60){$X'$}
\put(170,0){$W'$} \put(170,60){$B$} \put(112,65){\vector(1,0){53}}
\put(106,55){\vector(0,-1){41}} \put(175,55){\vector(0,-1){43}}
\put(114,58){\vector(1,-1){49}} \multiput(112,2.6)(5,0){11}{-}
\put(162,5){\vector(1,0){4}} \put(133,70){$f$} \put(180,30){$s$}
\put(92,30){$\pi$} \put(135,-8){$\varphi_{1}$}\put(136,40){$g$}
\end{picture}
\bigskip

We may write $K_{X'}=\pi^*(K_X)+E_{\pi}=_{\mathbb Q}M_1+Z_1,$ where
$M_1$ is the movable part of $|K_{X'}|$, $Z_1$ the fixed part and
$E_{\pi}$ an effective ${\mathbb Q}$-divisor which is a ${\mathbb
Q}$-sum of distinct exceptional divisors. By $K_{X'}- E_{\pi}$, we
mean $\pi^*(K_X)$. So, whenever we take the round up of
$m\pi^*(K_X)$, we always have $\roundup{m\pi^*(K_X)}\le mK_{X'}$ for
all positive numbers $m$. We may also write $\pi^*(K_X)=_{\mathbb Q}
M_1+E_1',$ where $E_1'=Z_1-E_{\pi}$ is actually an effective
${\mathbb Q}$-divisor.

If $\dim\varphi_{1}(X)=2$, a general fiber of $f$ is a smooth
projective curve of genus $\ge 2$. We say that $X$ is {\it
canonically fibred by curves}.

If $\dim\varphi_{1}(X)=1$, a general fiber $S$ of $f$ is a smooth
projective surface of general type. We say that $X$ is {\it
canonically fibred by surfaces} with invariants $(c_1^2(S_0),
p_g(S)),$ where $S_0$ is the minimal model of $S$. We may write
$M_1\equiv a_1S$ where $a_1\ge p_g(X)-1$.

{\it A generic irreducible element $S$ of} $|M_1|$ means either a
general member of $|M_1|$ whenever $\dim\varphi_{1}(X)\ge 2$ or,
otherwise, a general fiber of $f$.
\end{setup}

\begin{defn}
By abuse of terminology, we also define {\it a generic
irreducible element} $S'$ of an arbitrary linear system $|M'|$
on a general variety $V$ in a similar way. Assume that $|M'|$ is
movable. A generic irreducible element $S'$ is defined to be a
generic irreducible component in a general member of $|M'|$. (So
if $|M'|$ is composed with a pencil (i.e.
$\dim\Phi_{|M'|}(V)=1$), one has $M'\equiv tS'$ for some integer
$t\geq 1$. Clearly it may happen that sometimes $tS'\not\sim M'$
by our definition.)
\end{defn}

\section{\bf The key technical results}

In this section we will collect all the technical results needed
to prove our theorems. First we recall two lemmas about the
so-called Tankeev's principle which will be tacitly used
throughout our paper.

\begin{lem}\label{T} (\cite{T}, Lemma 2) Let $V$ be a
nonsingular projective variety. Let $D$ and $M$ are two divisors
on $V$. Assume that $|M|$ is base point free,
$\dim\Phi_{|M|}(V)\geq 2$ and $|D|\neq \emptyset$. Denote by $T$
a general member of $|M|$. If $\Phi_{|D+M|}$ is not birational,
then $\Phi_{|D+M|}|_T$ is not birational either.
\end{lem}

\begin{lem}\label{P12} (\cite{MPCPS}, 2.1) Let $V$ be a
nonsingular projective variety. Let $D$ and $M$ are two divisors
on $V$. Assume that $|M|$ is base point free,
$\dim\Phi_{|M|}(V)=1$ and $|D|\neq \emptyset$. Take the Stein
factorization of $\Phi_{|M|}$:
$$V\overset{f}\longrightarrow B\longrightarrow \bP^{h^0(V,M)-1}$$
where $f$ is a fibration onto the smooth curve $B$. Let $F$ be a
general fiber of $F$. Suppose we know (say by the vanishing
theorem) that $\Phi_{|D+M|}$ can separate different general
fibers of $f$ and that $\Phi_{|D+M|}|_F$ is birational for a
general $F$. Then $\Phi_{|D+M|}$ is birational.
\end{lem}

We will frequently use the following:

\begin{setup}\label{restrict}{\bf Fact.}
Let $|M|$ be a base point free linear system on any variety $V$
with $\dim\Phi_{|M|}(V)=1$. Denote by $S$ a generic irreducible
element of $|M|$. Then clearly ${\mathcal O}_S(M|_S)\cong
{\mathcal O}_S$ and ${\mathcal O}_S(S|_S)\cong {\mathcal O}_S$.
\end{setup}

{}From now on we present a technical, but key theorem which is a
generalized form of Theorem 2.6 in \cite{IJM} and will be
frequently applied in our proof.

\begin{setup}\label{J}{\bf Assumption}. We keep the same notation as
in \ref{1-map}. Assume that $X_0$ is a smooth model of $X$ and that
$J\le K_{X_0}$ is an effective divisor with $n_J:=h^0(X_0, \mathcal
{O}_{X_0}(J))\ge 2$. Recall that we have already a birational
morphism $\pi:X'\longrightarrow X$ in \ref{1-map}.

We may take further blow-ups to get a birational map $\pi_J:
X''\longrightarrow X$ such that the following three conditions are satisfied:
\par \vskip 1pc
\begin{itemize}
\item[$(i)_J$]
$X''$ is smooth and $X''$ dominates both $X'$ and $X_0$,
i.e. there are two birational morphisms $\pi':X''\longrightarrow X'$
and $\pi_0: X''\longrightarrow X_0$, $\pi_J:=\pi\circ \pi'$.
\item[$(ii)_J$]
The movable part $|M_J|$ of $|{\pi_0}^*(J)|$ is base point
free.
\item[$(iii)_J$]
The support of ${\pi_J}^*(K_1)$ is of simple normal
crossings.
\end{itemize}
\par \vskip 1pc
Denote by $g_J$ the composition $\phi_{|J|}\circ\pi_J$. So $g_J:
X''\longrightarrow W_J'\subseteq{\mathbb P}^{n_J-1}$ is a morphism.
Let $g_J: X''\overset{f_J}\longrightarrow
B_J\overset{s_J}\longrightarrow W_J'$ be the Stein factorization of
$g_J$. Since $M_J\le \pi'^*(M_1)\le \pi_J^*(K_X)$, we can write
$\pi_J^*(K_X)=_{\mathbb{Q}}M_J+E_J'$ where $E_J'$ is an effective
$\mathbb{Q}$-divisor. Set $d_J:=\dim (B_J)$. Denote by $S_J$ a
generic irreducible element of $|M_J|$. Then $S_J$ is a smooth
projective surface of general type. When $d_J=1$, one has $M_J\equiv
a_J S_J$ where $a_J\ge n_J-1$.
\par \vskip 1pc
(**) Whenever $p_g(X)=n_J$, we simply take $X'':=X'$ and adopt our
setting in \ref{1-map}. So $M_J$, $\pi_J$, $g_J$, $f_J$, $s_J$,
$S_J$, $B_J$, $W_J'$ and $E_J'$ are respectively $M_1$, $\pi$, $g$,
$f$, $s$, $S$, $B$, $W'$ and $E_1'$ just as in \ref{1-map}.
\end{setup}

\begin{setup}\label{main-notation}{\bf Key notations.} We place
ourselves in the above situation, and consider $S_J$, a generic
irreducible element of $|M_J|$. We define $p$ to be 1 if
$d_J\geq 2$ and $a_J$ otherwise. We assume further that $S_J$ is
equipped with a movable linear system $|G|$ and that a generic
irreducible element $C$ of $|G|$ is smooth. We define
$\xi:=\pi_J^*(K_X)\cdot C$. We fix an integer $m>0$, and we
consider a linear system on $S_J$, $L_m$, given by
$L_m:=|K_{S_J}+\roundup{(m-1)\pi_J^*(K_X)-S_J-\frac{1}{p}E_J'}|_{S_J}|$.
Then the following theorem holds:
\end{setup}

\begin{thm}\label{Key} (1) If $L_m$ separates different generic
irreducible elements of $|G|$ (namely $\Phi_{L_m}(C_1)\neq
\Phi_{L_m}(C_2)$ where $C_1$, $C_2$ are different generic
irreducible elements of $|G|$) and $\beta$ is a rational number
such that $\pi_J^*(K_X)-\beta C$ is numerically equivalent to an
effective ${\mathbb Q}$-divisor, then $\varphi_m$ is birational
if one of the following conditions is satisfied where we set
$\alpha:=(m-1-\frac{1}{p}-\frac{1}{\beta})\xi$ and
$\alpha_0:=\roundup{\alpha}$:

$i$. $\alpha > 2;$

$ii$. $\alpha_0\geq 2$ and $C$ is non-hyperelliptic;

$iii$. $\alpha>0$, $C$ is non-hyperelliptic and $C$ is an even
divisor on $S_J$.
\medskip

\noindent (2) One has the inequality $\xi\geq
\frac{2g(C)-2+\alpha_0}{m}$ if one of the following conditions
is satisfied:

$iv$. $\alpha > 1;$

$v$. $\alpha>0$ and $C$ is an even divisor on $S_J$.
\end{thm}
\begin{proof}
We consider the sub-system
$$|K_{X''}+\roundup{(m-1)\pi_J^*(K_X)-\frac{1}{p}E_J'}|\subset |mK_{X''}|.$$
Let $S_J'$ and $S_J''$ be two different generic irreducible elements
of $|M_J|$. Clearly one has
$$K_{X''}+\roundup{(m-1)\pi_J^*(K_X)-\frac{1}{p}E_J'}\ge
K_{X''}+\roundup{(m-2)\pi_J^*(K_X)}\ge M_J.$$ So
$|K_{X''}+\roundup{(m-1)\pi_J^*(K_X)-\frac{1}{p}E_J'}|$ can separate
$S_J'$ and $S_J''$ if either $\dim (B_J)\ge 2$ or $\dim (B_J)=1$ and
$g(B_J)=0$. For the case $\dim (B_J)=1$ and $g(B_J)>0$, one has
$a_J\ge n_J\ge 2$. Thus $p\ge 2$. Since
$$(m-1)\pi_J^*(K_X)-\frac{2}{p}E_J'-S_J'-S_J''\equiv
(m-1-\frac{2}{p})\pi_J^*(K_X)$$ is nef and big, the Kawamata-Viehweg
vanishing theorem gives a surjective map:
\begin{eqnarray*}
&& H^0(X'',K_{X''}+\roundup{(m-1)\pi_J^*(K_X)-\frac{2}{p}E_J'})\\
&\longrightarrow& H^0(S_J',
K_{S_J'}+\roundup{(m-1)\pi_J^*(K_X)-\frac{2}{p}E_J'}|_{S_J'})\oplus\\
&&H^0(S_J',
K_{S_J''}+\roundup{(m-1)\pi_J^*(K_X)-\frac{2}{p}E_J'}|_{S_J''}).
\end{eqnarray*}
The last two groups are non-zero. Therefore,
$$|K_{X''}+\roundup{(m-1)\pi_J^*(K_X)-\frac{1}{p}E_J'}|$$
separates $S_J'$ and $S_J''$. By Lemma \ref{T} and Lemma
\ref{P12}, it suffices to prove that $|mK_{X''}||_S$ gives a
birational map.

Noting that $(m-1)\pi_J^*(K_X)-\frac{1}{p}E_J'-S_J$ is nef and big,
the vanishing theorem gives a surjective map
\begin{eqnarray*}
&&H^0(X'',K_{X''}+\roundup{(m-1)\pi_J^*(K_X)-\frac{1}{p}E_J'})\\
&\longrightarrow  &H^0(S_J,
K_{S_J}+\roundup{(m-1)\pi_J^*(K_X)-S_J-\frac{1}{p}E_J'}|_{S_J}).
\end{eqnarray*}
We are reduced to prove that
$|K_{S_J}+\roundup{(m-1)\pi_J^*(K_X)-S_J-\frac{1}{p}E_J'}|_{S_J}|$
gives a birational map.

Now consider a generic irreducible element $C\in |G|$. By our
assumption there is an effective ${\mathbb Q}$-divisor $H$ on
$S_J$ such that
$$\frac{1}{\beta}\pi_J^*(K_X)|_{S_J}\equiv C+H.$$
By the vanishing theorem, we have the surjective map
$$H^0(S_J, K_{S_J}+\roundup{((m-1)\pi_J^*(K_X)-S_J-\frac{1}{p}E_J')|_{S_J}-H})\longrightarrow
H^0(C, K_C + D),$$ where
$D:=\roundup{((m-1)\pi_J^*(K_X)-S_J-\frac{1}{p}E_J')|_{S_J}-C-H}|_C$
is a divisor on $C$. Noting that
$$((m-1)\pi_J^*(K_X)-S_J-\frac{1}{p}E_J')|_{S_J}-C-H\equiv (m-1-\frac{1}{p}-\frac{1}{\beta})\pi_J^*(K_X)|_{S_J}$$
and that $C$ is nef on $S$, we have $\deg(D)\geq\alpha$ and thus
$\deg(D)\geq \alpha_0$. Whenever $C$ is non-hyperelliptic,
$m-1-\frac{1}{p}-\frac{1}{\beta}>0$ and $C$ is an even divisor
on $S$, $\deg(D)\geq 2$ automatically follows and thus $|K_C+D|$
gives a birational map.

Whenever $\deg(D)\ge 3$, then $|K_C+D|$ gives a birational map
and, by assumption the linear system $L_m$ separates different
irreducible elements of $|G|$, the Tankeev principle again says
that
$$|K_{S_J}+\roundup{((m-1)\pi_J^*(K_X)-S_J-\frac{1}{p}E_J')|_{S_J}-H}||_C$$
gives a birational map. Since
\begin{eqnarray*}
&&|K_{S_J}+\roundup{((m-1)\pi_J^*(K_X)-S_J-\frac{1}{p}E_J')|_{S_J}-H}|\\
&\subset
&|K_{S_J}+\roundup{(m-1)\pi_J^*(K_X)-S_J-\frac{1}{p}E_J'}|_{S_J}|,
\end{eqnarray*}
the latter linear system gives a birational map.
So $\varphi_m$ of $X$ is birational.

Whenever $\deg(D)\geq 2$, $|K_C+D|$ is base point free by the
curve theory. Denote by $|M_m|$ the movable part of $|mK_{X''}|$
and by $|N_m|$ the movable part of
$|K_{S_J}+\roundup{((m-1)\pi_J^*(K_X)-S_J-\frac{1}{p}E_J')|_{S_J}-H}|$.
By Lemma 2.7 of \cite{MPCPS}, we have
$$m\pi_J^*(K_X)|_{S_J}\ge N_m\ \ \text{and}\ \ (N_m\cdot C)_{S_J}\ge 2g(C)-2+\deg (D).$$
We are done.
\end{proof}

When applying Theorem \ref{Key}, another technical problem is
about how to find a suitable $\beta$ once we have chosen a
linear system $|G|$. The following lemma will be useful in this
context.
\begin{lem}\label{beta} Keep the same notation as in $\ref{J}$
and with $p$ as in $\ref{Key}$. Assume that $B_J=\mathbb{P}^1$. Let
$f_J:X''\longrightarrow \mathbb{P}^1$ be the induced fibration.
Denote by $F:=S_J$ a general fiber of $f_J$. Then one can find a
sequence of rational numbers $\{\beta_n\}$ with $\lim_{n\mapsto
+\infty} \beta_n = \frac{p}{p+1}$ such that
$\pi_J^*(K_X)|_F-\beta_n\sigma^*(K_{F_0}) \sim_{\bQ} N_n$ with an
effective ${\bQ}$-divisor $N_n$, where $\sigma:F\longrightarrow F_0$
is the blow down onto the smooth minimal model.
\end{lem}
\begin{proof} This is a generalized version of an established statement
of the first author (see Lemma 3.4 in \cite{Volume}).

One has $\mathcal {O}_{B_J}(p)\hookrightarrow {f_J}_*\omega_{X''}$
and therefore ${f_J}_*\omega_{X''/B_J}^{4p}\hookrightarrow
{f_J}_*\omega_{X''}^{4p+8}.$

For any positive integer $k$, denote by $M_k$ the movable part of
$|kK_{X''}|$. Note that ${f_J}_*\omega_{X''/B_J}^{4p}$ is generated
by global sections since it is semi-positive. So any local section
can be extended to a global one. On the other hand,
$|4p\sigma^*(K_{F_0})|$ is base point free and is exactly the
movable part of $|4pK_F|$ by Bombieri \cite{Bom} or Reider
\cite{Reider}. Applying Lemma 2.7 of \cite{MPCPS}, one has the
following, where $a_0:=4p+8$ and $b_0:=4p$:
$$a_0\pi_J^*(K_X)|_F\ge M_{4p+8}|_F\ge b_0\sigma^*(K_{F_0}).$$
This means that there is an effective $\mathbb{Q}$-divisor $E_0'$
such that
$$a_0\pi_J^*(K_X)|_F=_{\bQ} b_0\sigma^*(K_{F_0})+E_0'.$$
Thus $\pi_J^*(K_X)|_F =_{\bQ} \frac{p}{p+2}\sigma^*(K_{F_0})+E_0$
with $E_0=\frac{1}{a_0}E_0'$.

We first consider the case $p\ge 2$.

Assume that we have defined $a_n$ and $b_n$ such that the following
is satisfied with $l = n:$
$$a_{\ell}\pi_J^*(K_X)|_F \ge b_{\ell}\sigma^*(K_{F_0}).$$
We shall define $a_{n+1}$ and $b_{n+1}$ inductively such that the
above display is satisfied with $l = n+1$. One may assume from the
beginning that $a_n\pi_J^*(K_X)$ supports on a divisor with normal
crossings. Then the Kawamata-Viehweg vanishing theorem implies the
surjective map
$$H^0(K_{X''}+\roundup{a_n\pi_J^*(K_X)}+F)\longrightarrow H^0(F, K_F+
\roundup{a_n\pi_J^*(K_X)}|_F).$$ That means
\begin{eqnarray*}
|K_{X''}+\roundup{a_n\pi_J^*(K_X)}+F||_F&=&|K_F+\roundup{a_n\pi_J^*(K_X)}|_F|\\
&\supset& |K_F+b_n\sigma^*(K_{F_0})|\\
&\supset& |(b_n+1)\sigma^*(K_{F_0})|.
\end{eqnarray*}
Denote by $M_{a_n+1}'$ the movable part of $|(a_n+1)K_{X''}+F|$.
Applying Lemma 2.7 of \cite{MPCPS} again, one has
$M_{a_n+1}'|_F\ge (b_n+1)\sigma^*(K_{F_0}).$ \bigskip

{\bf Claim.} For any integer $t>0$,
$|K_{X''}+M_{a_n+t}'+F||_F\supset |(b_n+t+1)\sigma^*(K_{F_0})|$.
\begin{proof} Re-modifying our original $\pi_J$ such
that $|M_{a_n+1}'|$ is base point free. In particular,
$M_{a_n+1}'$ is nef. According to \cite{IJM}, $|mK_X|$ gives a
birational map whenever $m\ge 8$. Thus $M_{a_n+1}'$ is big. The
Kawamata-Viehweg vanishing theorem gives
\begin{eqnarray*}
|K_{X''}+M_{a_n+1}'+F||_F&=&|K_F+M_{a_n+1}'|_F|\\
&\supset& |K_F+(b_n+1)\sigma^*(K_{F_0})|\\
&\supset& |(b_n+2)\sigma^*(K_{F_0})|.
\end{eqnarray*}
Therefore $t=1$ case is true.

Assume we have proved that $|K_{X''}+M_{a_n+t-1}'+F||_F\supset
|(b_n+t)\sigma^*(K_{F_0})|$. Denote by $M_{a_n+t}'$ the movable
part of $|K_{X''}+M_{a_n+t-1}'+F|$. Then $M_{a_n+t}'|_F\supset
|(b_n+t)\sigma^*(K_{F_0})|$. Similarly we may assume that
$|M_{a_n+t}'|$ is base point free. The vanishing theorem gives
the surjective map:
\begin{eqnarray*}
|K_{X''}+M_{a_n+t}'+F||_F&=&|K_F+M_{a_n+t}'|_F|\\
&\supset& |K_F+(b_n+t)\sigma^*(K_{F_0})|\\
&\supset& |(b_n+t+1)\sigma^*(K_{F_0})|.
\end{eqnarray*}
\end{proof}

Take $t=p-1$. Noting that
$$|K_{X''}+M_{a_n+p-1}'+F|\subset |(a_n+p+1)K_{X''}|$$
and applying Lemma 2.7 of \cite{MPCPS} again, one has
$$a_{n+1}\pi_J^*(K_X)|_F\ge M_{a_n+p+1}|_F\ge M'_{a_n+p}|_F\ge
b_{n+1} \sigma^*(K_{F_0}).$$
Here we set $a_{n+1}:=a_n+p+1$ and
$b_{n+1}=b_n+p$. Note that
$a_n=n(p+1)+a_0$ and $b_n=np+b_0.$
Set $\beta_n = \frac{b_{n}}{a_{n}}.$
Then $\lim_{n\mapsto +\infty} \beta_n = \frac{p}{p+1}$.

The case $p=1$ can be proved similarly, but with a simpler
induction. We omit the details and leave it as an exercise.
\end{proof}
\section{\bf Proof of the main theorem, Part I}

We now begin the study of the birationality of $\varphi_4$. Set $ d
= d(X) = \dim\varphi_1(X)$. In the rest of the section, we shall
prove:
\begin{thm}\label{partI} Let $X$ be a minimal threefold of  general type with
$p_g(X)\ge 5$. Then $\varphi_4$ is {\bf not} birational if and only
if one of the following cases occurs, where $f:X'\longrightarrow B$
is as in Section $2$.
\begin{itemize}
\item[(i)]
$d(X) =2$, $g(C)=2$ and $\pi^*(K_X)\cdot C=1$ for a general
fiber $C$ of $f$.
\item[(ii)]
$d(X) = 1$ and $(K_{S_0}^2, p_g(S_0))=(1,2)$;
here $S_0$ is a smooth minimal model of a
general fiber $S$ of $f$.
\end{itemize}
\par
In both cases, $\varphi_4$ is generically finite of degree $2$ and
there is a family $\mathscr{C}$ (on $X$) of irreducible curves of geometric genus $2$ with
$K_X\cdot C_0=1$, where $C_0$ is a general member of $\mathscr{C}$.
\par
In the case $(ii)$, a general $C_0$ is the image on $X$
of a curve in the movable part of $|K_S|$ with $S$ a general fiber
of $f$.
\end{thm}
We prove the theorem according to the value of $d = d(X)$.

\begin{setup}{\bf The case $d = d(X) = 3$.}\label{d3} Assume $p_g(X)\ge 5$.
We shall show that $\varphi_4$ is birational.
\medskip

We set $J:=K_{X'}$ to run Theorem \ref{Key}. Then $\pi=\pi_J$.
According to the setting in Theorem \ref{Key}, we have $p=1$.
For a generic irreducible element $S$ of $|M_1|$, the linear
system $|M_1|_S|$ is not composed with a pencil of curves. Take
$G:=M_1|_S$. For all $m\ge 4$, it is clear that
$K_S+\roundup{(m-2)\pi^*(K_X)|_S}\ge G$. So Tankeev's principle
(Lemma \ref{T} and Lemma \ref{P12}) implies that
$|K_S+\roundup{(m-2)\pi^*(K_X)|_S}|$ separates different generic
irreducible elements of $|G|$. On the other hand, since
$\pi^*(K_X)|_S\ge M_1|_S\sim C$ where $C$ is an irreducible
curve on $S$. Take $\beta=1$ and the conditions of Theorem
\ref{Key}(1) is satisfied. One has $\pi^*(K_X)\cdot
S^2=\pi^*(K_X)|_S\cdot S|_S\ge ({M_1}|_S\cdot {M_1}|_S)=S^3$.
Denote by $\Lambda$ the subsystem $|S||_S\subset |S|_S|$. If
$\Phi_{\Lambda}$ is birational for a general $S$, so is
$\varphi_4$. Otherwise since $\Lambda$ gives a generically
finite map $S \rightarrow {\bP}^{p_g(X)-2}$ onto a
non-degenerate surface of degree $\ge 2$, one has $(S|_S)^2\ge
2(p_g(X)-3)\ge 4$. So $\xi=\pi^*(K_X)\cdot S^2\ge 4$.

If we take $m=4$, then
$$\alpha=(m-1-\frac{1}{p}-\frac{1}{\beta})\xi\ge \xi>2.$$
Theorem \ref{Key}(1) tells us that $\varphi_4$ is birational
onto its image. (This method, however, gives no information on
the birationality of $\varphi_3$.)
\end{setup}

\begin{setup}\label{d2}{\bf The case $d = d(X) =2$.} Assume $p_g(X)\ge 5$.
We shall show by the arguments up to \ref{6/5} that $\varphi_4$ is
birational unless the case of Theorem \ref{partI} (i) occurs. Set
$J:=K_{X'}$. Let $S$ be a generic irreducible element of $|M_1|$.
\medskip

By our definition in \ref{main-notation}, one has $p=1$. Note
that $\pi^*(K_X)|_S\ge {M_1}|_S\equiv a_2C$, where $C$ is a
general fiber of $f:X'\longrightarrow B$ and $a_2\geq h^0(S,
M_1|_S)-1\geq p_g(X)-2$. Take $G:={M_1}|_S$. Then $C$, as a
generic irreducible element of $|G|$, is smooth. So we can take
$\beta:=a_2\ge 3$. Pick up two different general fibers $C_1$
and $C_2$ of $f$. Since we have
$$\pi^*(K_X)|_S\equiv a_2C+E_1'|_S,$$
for all $m\ge 4$,
$$(m-2)\pi^*(K_X)|_S-C_1-C_2-\frac{2}{a_2}E_1'|_S\equiv
(m-2-\frac{2}{a_2})\pi^*(K_X)|_S$$ is nef and big. So the
Kawamata-Viehweg vanishing theorem (\cite{KV, V}) gives a surjective
map
\begin{eqnarray*}
& &H^0(S,
K_S+\roundup{(m-2)\pi^*(K_X)|_S-\frac{2}{a_2}E_1'|_S})\\
&\longrightarrow &H^0(C_1, K_{C_1}+D_1)\oplus H^0(C_2, K_{C_2}+D_2),
\end{eqnarray*}
where
$D_i:=\roundup{(m-2)\pi^*(K_X)|_S-\frac{2}{a_2}E_1'|_S}|_{C_i}$
for $i=1,2$, $H^0(C_i, K_{C_i}+D_i)\ne 0$ because $\deg(D_i)\ge
(m-2-\frac{2}{a_2})(\pi^*(K_X)|_S\cdot C_i)>0$ for $i=1,2$. This
means that
$|K_S+\roundup{(m-2)\pi^*(K_X)|_S-\frac{2}{a_2}E_1'|_S}|$
separates $C_1$ and $C_2$. So does
$|K_S+\roundup{(m-2)\pi^*(K_X)|_S}|$. Thus conditions in Theorem
\ref{Key}(1) are always satisfied. Now we may run Theorem
\ref{Key}.

Take a sufficiently large integer $m$ such that $\alpha>1$, one has
then the inequality
$$m\xi\ge 2g(C)-2+(m-1-\frac{1}{p}-\frac{1}{\beta})\xi,$$
which gives
$$\xi\ge \frac{3(2g(C)-2)}{7}.\eqno (4.1)$$

Take $m=4$. Then
$$\alpha=(2-\frac{1}{\beta})\xi\ge \begin{cases} \frac{10}{7}>1,\ \ &
g(C)=2;\\
\frac{20}{7}>2,\ \ & g(C)\ge 3.\end{cases}$$

So Theorem \ref{Key} says that $\xi\ge 1$ whenever $g(C)=2$; that
$\xi\ge \frac{7}{4}$, $\varphi_4$ is birational whenever $g(C)\ge
3$. We will study the case $g(C)=2$ which is slightly complicated.

\begin{claim}\label{g2xi>1} (1) If $g(C)=2$, $p_g(X)\ge 6$ and $\xi>1$, then $\xi\ge
\frac{5}{4}$ and $\varphi_4$ is birational onto its image.

(2) If $g(C)=2$, $p_g(X)=5$ and $\xi>1$, then $\xi\ge \frac{6}{5}$.
Whenever $\xi> \frac{6}{5}$, $\varphi_4$ is birational onto its
image (which also implies $\xi\ge \frac{5}{4}$).
\end{claim}
\begin{proof} (1) We have $p=1$ and $\beta\ge 4$ by assumption.
We can always find a positive integer $l_0>4$ such that $\xi\ge
\frac{l_0+1}{l_0}$. Take $m_0=l_0-1\ge 4$. We hope to run Theorem
\ref{Key}. We have:
\begin{eqnarray*}
\alpha-(l_0-3)&=&(l_0-2-\frac{1}{p}-\frac{1}{\beta})\xi-(l_0-3)\\
&\ge& (l_0-3-\frac{1}{4})\frac{l_0+1}{l_0}-(l_0-3)\\
&=& \frac{1}{4}(3-\frac{13}{l_0})>0.
\end{eqnarray*}
Thus one has $\alpha_0\ge l_0-2$. Theorem \ref{Key} gives
$\xi\ge\frac{l_0}{l_0-1}$. One may proceed by induction as long as
$l_0-1 \ge 4$. So Theorem \ref{Key} simply gives $\xi\ge
\frac{5}{4}$.

Take $m_1=4$. We get $\alpha=(4-2-\frac{1}{\beta})\xi\ge
\frac{7}{4}\cdot \frac{5}{4}>2.$ Therefore $\varphi_4$ is birational
onto its image.

(2) The same argument as in (1) shows that if $\beta\geq 3$ then
$\alpha-(l_0-3)\geq \frac{2l_0-10}{3l_0}$ and hence
$\alpha_0\geq l_0-2$ if $l_0>5$. If we take $l_0=6$ and
$m_0=l_0-1=5$. Then $\alpha_0\geq 4$. Theorem \ref{Key}(2) gives
$\xi\geq\frac{6}{5}$.

Take $m_2=4$. We get $\alpha=(4-2-\frac{1}{\beta})\xi\ge
\frac{5}{3}\cdot \frac{6}{5}=2.$ Therefore, $\varphi_4$ is
birational onto its image by Theorem \ref{Key} whenever
$\xi>\frac{6}{5}$. Furthermore Theorem \ref{Key}(2) gives
$\xi\geq \frac{5}{4}$ whenever $\xi>\frac{6}{5}$.
\end{proof}

\begin{claim}\label{6/5} The situation with $d=2$, $g(C)=2$, $p_g(X)=5$ and
$\xi=\frac{6}{5}$ does not occur.
\end{claim}
\begin{proof} Assume $d=2$, $g(C)=2$, $p_g(X)=5$ and
$\xi=\frac{6}{5}$. We want to deduce a contradiction.

Consider a general member $S\in |M_1|$. We have seen that
$M_1|_S\equiv a_2C$ where $a_2\ge p_g(X)-2\ge 3$. If $a_2>3$, then
we may take $\beta=4$ as in the proof of \ref{g2xi>1}. Then Theorem
\ref{Key} will give $\xi>\frac{6}{5}$, a contradiction. So we may
assume $a_2=3$. Consider the induced fibration $f:X'\longrightarrow
B$ where $B$ is a normal surface given by the Stein factorization;
see the notation in Section 2. We can write $S=f^*(H_B)$ for some
ample divisor $H_B$ on $B$ and $H_B=s^*(H')$ for a hyperplane $H'$
on $W'\subset \mathbb{P}^4$. Then $3=a_2=H_B^2=\deg(s) \deg (W')$.
Since $W'$ is non-degenerate, $\deg (W')\ge 3$. This means that
$\deg (s)=1$ and $s$ is a finite morphism of degree 1. Hence $B$ is
the normalization of $W'$. According to Del Pezzo, or Nagata (see
Theorem 7 in \cite{Na}), or Reid (see Exercise 19 of Chapter 2, page
30 in \cite{Park}), $W'$ is either a cone over a smooth rational
base curve of degree 3 in $\mathbb{P}^3$, or the ruled surface
$\mathbb{F}_e$. In particular, $W'$ is a normal rational surface and
$W'$ has at worst a single singularity. We know that a minimal
resolution $\bar{W}$ of $W'$ is the ruled surface $\mathbb{F}_e$
with $e\ge 0$. Denote by $\rho: \bar{W}\longrightarrow W'$ the
minimal resolution. Set $\bar{H}=\rho^*(H')$. Then $\bar{H}^2=
(H')^2 = \deg(W')=3$. So it is clear that $\bar{W}$ can not be
$\mathbb{F}_0$. Thus we have $\bar{W}=\mathbb{F}_e$ with $e>0$.

If necessary we can remodify our $X'$ (by further blowing ups) so
that we get a fibration $f:X'\longrightarrow B$ with $B$ smooth so
that $g = s \circ f$ in notation of Section 2; further, $s = \rho
\circ \tau$, where $\tau:B\longrightarrow \mathbb{F}_e$ and
$\rho:\mathbb{F}_e \longrightarrow W'$. We have $H_B=\tau^*\bar{H}$.
Now we can perform the computation on $\mathbb{F}_e$ with $e>0$.
Noting that $\bar{H}$ is nef and big on $\mathbb{F}_e$, we can write
$$\bar{H}\sim \mu G_0+nT$$
where $G_0$ is the unique section with $G_0^2=-e$, $\mu$ and $n$ are
integers and $T$ is the general fiber of the ruling on
$\mathbb{F}_e$. The property of $\bar{H}$ being nef and big implies
$\mu>0$ and $n\ge \mu e$. If $e=1$, the equality $(\mu G_0+nT)^2=3$
implies $\mu=1$ and $n=2$. If $e>1$, clearly $n\ge 2\mu\ge 2$. Now
let $\alpha_e:\mathbb{F}_e\longrightarrow \mathbb{P}^1$ be the
ruling, whose fibers are all smooth rational curves. Set
$f_0:=\alpha_e\circ\tau\circ f: X'\longrightarrow \mathbb{P}^1$,
which is a fibration with connected fibers. Denote by $F$ a general
fiber of $f_0$. We have
$$M_1\sim f^*(H_B)=f^*\tau^*(\bar{H})\ge nF$$
with $n\ge 2$.

Note that a general fiber $F$ of $f_0$ is fibred by curves $C$ (also
as fibers of $f$) with $g(C)=2$. In our situation, we may take $p=2$
and consider $J:=2F$ on $X'$. This fits the setting for Lemma
\ref{beta}. Note that we actually have $f_J=f_0$ and $\pi_J=\pi$
under this situation. So Lemma \ref{beta} says that
$\pi^*(K_X)|_F-\beta_0\sigma^*(K_{F_0})$ is pseudo-effective for
$\beta_0\mapsto\frac{2}{3}$. One has
$$\frac{6}{5}=\pi^*(K_X)\cdot C=(\pi^*(K_X)|_F\cdot C)_F\ge
\frac{2}{3}\sigma^*(K_{F_0})\cdot C.$$ This implies that
$\sigma^*(K_{F_0})\cdot C=1$.  Note that the uniqueness of the
Zariski decomposition says that the nef divisor $\pi^*(K_X)|_F$ can be added
some effective ${\mathbb Q}$-divisors to become
the maximal nef part $\sigma^*(K_{F_0})$ in $K_F$. Thus $\sigma^*(K_{F_0})-\pi^*(K_X)|_F$ is
pseudo-effective. One has $\frac{6}{5}=(\pi^*(K_X)|_F\cdot
C)_F\le \sigma^*(K_{F_0})\cdot C=1$, which is absurd. We are
done.
\end{proof}

\begin{prop}\label{g2xi1} If $g(C)=2$, $p_g(X)\ge 5$ and $\xi=1$, then $\varphi_4$
is generically finite of degree 2.
\end{prop}
\begin{proof}
Recall that we have $K_{X'}=\pi^*(K_X)+E_{\pi}$. On $X$ we set
$Z:=\pi_*(Z_1)$ and $N:=\pi_*(M_1)$. Clearly $K_X\sim N+Z$. Then
there is an effective $\mathbb{Q}$-divisor $E_1$, which is
supported by some exceptional divisors, such that
$\pi^*(N)=M_1+E_1$. Therefore $E_1'=\pi^*(Z)+E_1$. For a general
member $S$ of $|M_1|$, we have
$K_{X'}|_S=\pi^*(K_X)|_S+E_{\pi}|_S=({M_1}|_S+E_1'|_S)+E_{\pi}|_S$.
So ${E_1'}|_S=\pi^*(Z)|_S+{E_1}|_S$. One knows that $E_{\pi}$ is
composed of all those exceptional divisors of $\pi$. Thus
$E_1\leq E_{\pi}$ and $E_1|_S\leq E_{\pi}|_S$.

We will now need further assumptions on the map $\pi$. we may
take the $\pi$ to be the composition
$X'\overset{\pi_2}\longrightarrow
X_2\overset{\pi_1}\longrightarrow
X_1\overset{\pi_0}\longrightarrow X$ where $\pi_0$ is the
resolution of the indeterminancy of $\varphi_1$, $\pi_1$ is the
resolution of those isolated singularities on $X_1$ which are
away from all exceptional locus of $\pi_0$ and $\pi_2$ is the
minimal further modification such that $\pi^*(K_1)$ has the
support of simple normal crossings (recall here that $K_1\sim
K_X$ is a fixed Weil divisor as in \ref{1-map}). Set
$\pi_3:=\pi_0\circ \pi_1$. By abuse of notations we will have a
set of divisors for $\pi_3$ similar to that for $\pi$. For example we
may write $K_{X_2}=\pi_{3}^*(K_X)+E_{\pi_3}$ where $E_{\pi_3}$
is an effective ${\mathbb Q}$-divisor. The movable part
$|M_{\pi_3}|$ of $|K_{X_2}|$ is already base point free. Write
$\pi_3^*(N)=M_{\pi_3}+E_{1, \pi_3}$ and
$\pi_3^*(K_X)=M_{\pi_3}+E_{1,\pi_3}'$ where $E_{1, \pi_3}$ and
$E_{1, \pi_3}'$ are both effective ${\mathbb Q}$-divisors.
Clearly $E_{1, \pi_3}'=\pi_3^*(Z)+E_{1, \pi_3}$. By the
definition of $\pi_3$, $E_{\pi_3}$ is the union of two parts
$E_{\pi_3}'+E_{\pi_3}''$ where $E_{\pi_3}'$ consists all those
components over the indeterminancy of $\varphi_1$ and
$E_{\pi_3}''$ is totally disjoint from $E_{\pi_3}'$. Denote by
$S_{\pi_3}$ a general member of $|M_{\pi_3}|$. Then
$|M_{\pi_3}|_{S_{\pi_3}}|$ is a free pencil of genus 2 with a
general member $C_{\pi_3}$. As we have seen
$\text{Supp}(E_{\pi_3}''|_{S_{\pi_3}})=0$ and so
$\text{Supp}(E_{\pi_3}|_{S_{\pi_3}})=
\text{Supp}(E_{1,\pi_3}|_{S_{\pi_3}})$. Now we see
$1=\pi^*(K_X)\cdot C=\pi_3^*(K_X)\cdot
{\pi_2}_*(C)=\pi_3^*(K_X)\cdot C_{\pi_3}.$ Since $2 =
\deg(K_{C_{\pi_3}})= (\pi_3^* K_X + E_{\pi_3}) \cdot C_{\pi_3}$
and $\pi_3^*(K_X)|_{S_{\pi_3}}\cdot C_{\pi_3}=1$, we get
$E_{\pi_3}|_{S_{\pi_3}}\cdot C_{\pi_3}=E_{\pi_3}\cdot
C_{\pi_3}=1$. Therefore $E_{1,\pi_3}\cdot C_{\pi_3}>0$. Noting
that $\pi_2^*(E_{1,\pi_3})\leq E_{1}$ one sees
$$E_{1}|_S\cdot C\geq \pi_2^*(E_{1,\pi_3})|_S\cdot
C=E_{1,\pi_3}|_{S_{\pi_3}}\cdot C_{\pi_3}>0 \eqno{(4.2)}$$ which
is what we want to show in this paragraph.

Since $d(X) > 1$, the linear system $|4K_{X'}|$ (and also
$|K_{X'}|$) separates different irreducible elements of $|M_1|$.
Therefore, $\Phi_4$ is birational if and only if $\Phi_4|_S$ is
birational for a general $S$. Now on a general surface $S$, we
have a pencil $|M_1|_S|$ and $\Phi_4|_S$ separates different
generic irreducible elements of $|M_1|_S|$ (see \ref{d2} and the
paragraph below). So similarly $\Phi_4|_S$ is birational if and
only if $(\Phi_4|_S)|_C$ is birational. We will show that
$(\Phi_4|_S)|_C=\Phi_{|2K_C|}$. The latter is, however, not
birational.

Denote by $M_4$ the movable part of $|4K_{X'}|$. Then, for a general
$S$ in $|M_1|$, $M_4|_S\le 4\pi^*(K_X)|_S$ and $M_4|_S\cdot C\le
4\pi^*(K_X)|_S\cdot C=4$. By Kawamata-Viehweg vanishing, we have a
surjective map:
$$H^0(X', K_{X'}+\roundup{2\pi^*(K_X)}+S)\longrightarrow H^0(S,
K_S+\roundup{2\pi^*(K_X)}|_S).$$ Since
$$2\pi^*(K_X)|_S-C-\frac{1}{a_2}E_1'|_S\equiv
(2-\frac{1}{a_2})\pi^*(K_X)|_S$$ is nef and big, the vanishing
theorem gives a surjective map:
$$H^0(S,K_S+\roundup{2\pi^*(K_X)|_S-\frac{1}{a_2}E_1'|_S})\longrightarrow
H^0(C, K_C+D),$$ where
$$D:=\roundup{2\pi^*(K_X)|_S-\frac{1}{a_2}E_1'|_S}|_C=\roundup{(2-\frac{1}{a_2})E_1'|_S}|_C$$
and $\deg(D)\ge (2-\frac{1}{a_2})\ge 2-\frac{1}{3}>1$, noting that
$E_1' \cdot C = \pi^* K_X \cdot C = \xi = 1$. So $|K_C+D|$ is base
point free. Denote by $M_4'$, $N_4'$ the movable parts of
$|K_{X'}+\roundup{2\pi^*(K_X)}+S|$,
$|K_S+\roundup{2\pi^*(K_X)|_S-\frac{1}{a_2}E_1'|_S}|$ respectively.
Then, by Lemma 2.7 of \cite{MPCPS}, one has
$$M_4|_C:=(M_4|_S)|_C\ge (M_4'|_S)|_C\ge N_4'|_C\ge K_C+D.$$
So $4=4\pi^*(K_X)|_S\cdot C \ge M_4\cdot C=\deg(K_C+D)\ge 4$. This means
$M_4|_C\sim K_C+D$ and $\deg(D) = 2$. On the other hand, we have shown
$|M_4||_C\supset |K_C+D|$. Thus $\Phi_4|_C=\Phi_{|K_C+D|}$.
Since $\deg(\Phi_{|K_C|}) = 2$, we have
$\deg(\Phi_4) \le 2$. So $\Phi_4$ is either
birational or of degree two.

We have:
$$K_C\sim (K_{X'}|_S + S|_S)|_C=(\pi^*(Z)|_S)|_C+(E_1|_S)|_C+(E_{\pi}|_S)|_C.
\eqno {(4.2)}$$ Since $2 = \deg(K_C)= (\pi^* K_X + E_{\pi})
\cdot C$ and $E_1'|_S\cdot C=\pi^*(K_X)|_S\cdot C=1$, we get
$E_{\pi}|_S\cdot C=E_{\pi}\cdot C=1$. So
\text{Supp}($(E_{\pi}|_S)|_C)$ is in either of the following
situations:

{\bf Case 1.} a single point $P$;

{\bf Case 2.} two points $P+Q$, where $P$, $Q$ are different
points on $C$.

We consider {\bf Case 1} and {\bf Case 2} separately and note
that $E_1'|_S=\pi^*(Z)|_S+E_1|_S$ and \text{Supp}($E_1|_S$)
$\subset$ \text{Supp}($E_{\pi}|_S$).

Suppose we are in {\bf Case} 1. Then $(E_{\pi}|_S)|_C=P$. If
\text{Supp}($(\pi^*(Z)|_S)|_C+E_1|_S$) contains a point other
than $P$ (say a point R), then
$\pi^*(Z)|_S)|_C+(E_1|_S)|_C+(E_{\pi}|_S)|_C=P+R$ and $R$ is not
contained in \text{Supp}$(E_1|_S)$ otherwise $R$ is in
\text{Supp}($E_{\pi}|_S$), a contradiction. Thus $R\leq
\pi^*(Z)|_S)|_C$ as an integral part because
$\pi^*(Z)|_S)|_C+(E_1|_S)|_C+(E_{\pi}|_S)|_C$ is an integral
divisor. This says that $\deg(\pi^*(Z)|_S)|_C+(E_1|_S)|_C)>1$
because $\deg((E_1|_S)|_C)>0$ by the relation (4.2), which is
absurd. Thus $(\pi^*(Z)|_S)|_C+(E_1|_S)|_C=P$ and $K_C\sim 2P$.
In this case, we have $D=2P$ and $\varphi_4|_C=\Phi_{|2K_C|}$ is
not birational. So $\varphi_4$ is not birational onto its image.

Suppose we are in {\bf Case} 2. The right hand side of (4.2)
must be $P+Q$ and $K_C\sim P+Q$. We also know that
Supp$((\pi^*(Z)|_S)|_C+(E_1|_S)|_C)=P+Q$. So $D=P+Q$. And thus
$\varphi_4|_C=\Phi_{|2K_C|}$ is not birational.
%
%
%
%
\end{proof}
\end{setup}

{}From now on, we study the case $d = d(X) =1$. We shall show that
$\varphi_4$ is birational unless the case of Theorem \ref{partI}
(ii) occurs. Let $b$ be the genus of $B = f(X')$. Let $S$ be a
general fibre of $f : X' \rightarrow B$ and $\sigma : S \rightarrow
S_0$ the smooth blow down to a minimal model. {}From now on within
this section, we always set $J:=K_{X'}$ to run Theorem \ref{Key}.

\begin{lem}\label{b>0} If $d=1$ and $b>0$, then
$\pi^*(K_X)|_S\sim\sigma^*(K_{S_0})$.
\end{lem}
\begin{proof} We use the idea of Lemma 14 in Kawamata's paper \cite{KA}.
By Shokurov's theorem in \cite{Sho}\footnote{C. D. Hacon and J.
M$^{\rm c}$Kernan (see \cite{HsM}) have recently extended V. V.
Shokurov's result to any dimension and without assuming the MMP.},
each fiber of $\pi:X'\longrightarrow X$ is rationally chain
connected. Therefore, $f(\pi^{-1}(x))$ is a point for all $x\in X$.
Considering the image $G\subset (X \times B)$ of $X'$ via the
morphism $(\pi\times f)\circ \triangle_{X'}$ where $\triangle_{X'}$
is the diagonal map $X'\longrightarrow X'\times X'$, one knows that
$G$ is a projective variety. Let $g_1:G\longrightarrow X$ and
$g_2:G\longrightarrow B$ be two projections. Since $g_1$ is a
projective morphism and even a bijective map, $g_1$ must be both a
finite morphism of degree 1 and a birational morphism. Since $X$ is
normal, $g_1$ must be an isomorphism. So $f$ factors as $f_1 \circ
\pi$ where $f_1:=g_2\circ g_1^{-1} : X \rightarrow B$ is a well
defined morphisms. In particular, a general fiber $S_0$ of $f_1$
must be smooth minimal. So it is clear that
$\pi^*(K_X)|_S\sim\sigma^*(K_{S_0})$ where $\sigma$ is nothing but
$\pi|_S$.
\end{proof}

\begin{setup}\label{d1b>0} {\bf The case $d=1$, $b>0$.} Assume
$p_g(X)\ge 2$ and $S$ is not of type $(c_1^2, p_g) = (1,2)$.
We shall show that $\varphi_4$ is birational.
\medskip

One has $M_1\equiv a_1S$ where $a_1\ge p_g(X)\ge 2$. Clearly,
$|4K_{X'}|$ separates different generic irreducible elements of
$|M_1|$. In fact, if $S_1$ and $S_2$ are two different fibers of
$f$, Kawamata-Viehweg vanishing gives a surjective map
\begin{eqnarray*}
& &H^0(X',
K_{X'}+\roundup{3\pi^*(K_X)-\frac{2}{a_1}E_1'})\\
&\longrightarrow & H^0(S_1,
K_{S_1}+\roundup{3\pi^*(K_X)-\frac{2}{a_1}E_1'}|_{S_1})\oplus\\
& &H^0(S_2, K_{S_2}+\roundup{3\pi^*(K_X)-\frac{2}{a_1}E_1'}|_{S_2}).
\end{eqnarray*}
One has
\begin{eqnarray*}
&&K_{S_i}+\roundup{3\pi^*(K_X)-\frac{2}{a_1}E_1'}|_{S_i}\\
&\ge&
K_{S_i}+\roundup{2\pi^*(K_X)|_{S_i}+(1-\frac{2}{a_1})E_1'|_{S_i}}\\
&\ge& K_{S_i}+2\sigma_i^*(K_{S_{i}})>0
\end{eqnarray*}
which means that $|K_{X'}+\roundup{3\pi^*(K_X)-\frac{2}{a_1}E_1'}|$
separates $S_1$ and $S_2$. So does $|4K_{X'}|$.

We study another subsystem
$$|K_{X'}+\roundup{3\pi^*(K_X)-\frac{1}{a_1}E_1'}|\subset
|4K_{X'}|.$$  Pick up a general fiber $S$ of $f: X'\longrightarrow
B$. Since $3\pi^*(K_X)-\frac{1}{a_1}E_1'-S\equiv
(3-\frac{1}{a_1})\pi^*(K_X)$ is nef and big, the Kawamata-Viehweg
vanishing theorem gives a surjective map
\begin{eqnarray*}
& &H^0(X',
K_{X'}+\roundup{3\pi^*(K_X)-\frac{1}{a_1}E_1'})\\
&\longrightarrow &
H^0(S, K_S+\roundup{3\pi^*(K_X)-\frac{1}{a_1}E_1'}|_S)\\
&\supset& H^0(S, K_S+\roundup{(3\pi^*(K_X)-\frac{1}{a_1}E_1')|_S})\\
&=& H^0(S,
K_S+2\sigma^*(K_{S_0})+\roundup{(1-\frac{1}{a_1})E_1'|_S}).
\end{eqnarray*}
Thus, by Tankeev's principle, it suffices to show the birationality
of
$\Phi_{|K_S+2\sigma^*(K_{S_0})+\roundup{(1-\frac{1}{a_1})E_1'|_S}|}.$

If $(K_{S_0}^2, p_g(S_0))\ne (2,3)$, the statement is clear by
Bombieri (\cite{Bom}). Otherwise it is a corollary of Lemma 1.3 in
\cite{JPAA}. (In fact, this is an easy exercise by Kawamata-Viehweg
vanishing, just noting the fact (see page 227 in \cite{BPV}) that
$|K_{S_0}|$ has no base points and that $|K_{S_0}|$ gives a
generically finite morphism of degree 2.)
\end{setup}

\begin{setup}\label{d1b=0}
{\bf Claim.} Consider the case $d=1$, $b=0$. Assume that
$p_g(X)\ge 5$ and that a general fiber $S$ of the induced fibration
$f: X'\longrightarrow B$ is not of type (1,2). Then
$\varphi_4$ is birational.

Since $p_g(X)>0$, we have $p_g(S)>0$ for a general fiber $S$. Thus
by an established theorem (see Bombieri \cite{Bom}, Reider
\cite{Reider}, Catanese-Ciliberto \cite{C-C}, and P. Francia
\cite{Francia} or directly refer to Theorem 3.1 in the survey
article by Ciliberto \cite{Ci}), $|2K_{S_0}|$ is always base point
free. We classify $S$ into two types separately in order to organize
our proof ((1,2) type surfaces excluded):

(A) $(K_{S_0}^2, p_g(S))\ne (2,3)$;

(B) $(K_{S_0}^2, p_g(S))= (2,3)$.

To prove the claim, we consider these two case separately.

{\bf Case (A)}. We want to verify the conditions in Theorem
\ref{Key}. We take $G:=2\sigma^*(K_{S_0})$. Since $|G|$ is base
point free, a generic irreducible element $C$ of $|G|$ is smooth
and $|G|$ is not composed with a pencil. By the assumption, we
have an inclusion ${\mathcal{O}}(4)\hookrightarrow
f_*\omega_{X'}$. So there is another inclusion
$$f_*\omega^2_{X'/\mathbb{P}^1}\hookrightarrow
f_*\omega^3_{X'}.\eqno {(4.3)}$$ The sheaf on the left is
semi-positive and is clearly generated by global sections. On the
other hand, $3\pi^*(K_X)\ge M_3$ where $M_3$ is the movable part of
$|3K_{X'}|$. Noting that any local sections along a general fiber
$S$ can be extended to a global section, we clearly have
$3\pi^*(K_X)|_S\ge 2\sigma^*(K_{S_0})\sim G$. So we have
$$K_S+\roundup{3\pi^*(K_X)-\frac{1}{a_1}E_1'}|_S\ge
\roundup{E_{\pi}|_S+(1-\frac{1}{a_1})E_1'|_S+3\pi^*(K_X)|_S}\ge
G.$$ So one of the conditions in Theorem \ref{Key}(1) is
satisfied.

By Lemma \ref{beta}, we may take a $\beta$ (arbitrarily close to
$\frac{2}{5}$) such that $\pi^*(K_X)|_S-\beta C$ is
pseudo-effective. Noting that $G$ is an even divisor and that $C$ is
non-hyperelliptic (by the birationality of $\Phi_{|3K_S|}$), we only
have to verify $\alpha>0$ in order to apply Theorem \ref{Key}. We
have $p=4$. Take $m=4$. Then
$$\alpha=(4-1-\frac{1}{4}-\frac{1}{\beta})\xi>0.$$
So $\varphi_4$ is birational by Theorem \ref{Key}.

{\bf Case (B)}. We take $G=\sigma^*(K_{S_0})$. Since $|G|$ is
base point free, a generic irreducible element $C$ of $|G|$ is a
smooth curve of genus 3. Since $|2\sigma^*(K_{S_0})|$ is base
point free, the similar argument in Case (A) shows that
conditions in Theorem \ref{Key}(1) are satisfied.

We have $p=4$. Again by Lemma \ref{beta}, we may take a $\beta$
(arbitrarily near $\frac{4}{5}$) such that $\pi^*(K_X)|_S-\beta C$
is  pseudo effective. We have $\xi\ge \beta C^2\mapsto \frac{8}{5}$.
In fact, taking the limit, one has $\xi\ge \frac{8}{5}$. Take $m=4$.
Then $\alpha=(4-1-\frac{1}{p}-\frac{1}{\beta})\xi\ge\frac{12}{5}>2$.
Theorem \ref{Key} tells us that $\varphi_4$ is birational onto its
image. This proves the claim.
\end{setup}

\begin{setup}\label{d1(1,2)}{\bf The case $d=1$, $(K_{S_0}^2,
p_g(S_0))=(1,2)$.} Assume $p_g(X)\ge 3$.
We shall show that the case of Theorem \ref{partI} (ii) occurs.
\medskip

Since the 4-canonical map of $S$ is not birational by Bombieri (see
\cite{Bom} or \cite{BPV}), it is clear that $\varphi_4$ of $X$ is
not birational either. It is however easy to see by inclusion (4.3)
that $\varphi_4$ is at worst generically finite of degree 2 since
$\Phi_{|2K_S|}$ is generically finite of degree 2.

We now show the existence of the family $\mathscr{C}$ of curves of
genus 2 on $X$ with the property $K_X\cdot C_0=1$.

We may consider the relative canonical map $\Psi: X'\dashrightarrow
\mathbb{P}(f_*\omega^{\vee}_{X'/B})$ over $B$. By taking further
birational modifications we may assume that $\Psi$ is a morphism
over $B$. So we have the following commutative diagram:
\medskip

\begin{picture}(50,80) \put(100,0){$X$} \put(100,60){$X'$}
\put(170,0){$B$} \put(170,60){$\mathbb{P}(f_*\omega^{\vee}_{X'/B})$}
\put(112,65){\vector(1,0){53}} \put(106,55){\vector(0,-1){41}}
\put(175,55){\vector(0,-1){43}} \put(114,58){\vector(1,-1){49}}
\multiput(112,2.6)(5,0){11}{-} \put(162,5){\vector(1,0){4}}
\put(133,70){$\Psi$} \put(180,30){$p$} \put(92,30){$\pi$}
\put(135,-8){$\varphi_{1}$}\put(136,40){$f$}
\end{picture}
\bigskip

Since $\Psi|_S=\Phi_{|K_S|}$, we see that a  general fiber of $\Psi$
is a genus 2 curve which is nothing but the movable part of $|K_S|$.
By \cite{BPV}, one knows that $|K_{S_0}|$ has one base point and its
generic irreducible element is a curve of genus 2. So on $S$ we take
$G$ to be the movable part of $|\sigma^*(K_{S_0})|$. Pick up a
generic irreducible element $C$ of $|G|$, which is also a smooth
fiber of $\Psi$. We will apply Theorem \ref{Key} to estimate
$\xi=\pi^*(K_X)|_S\cdot C$. Note that $\xi\le \sigma^*(K_{S_0})
\cdot C = K_{S_0}^2 = 1$ because $\pi^*(K_X)|_S$ is nef and
$\sigma^*(K_{S_0})$ is the nef part of the Zariski decomposition of
$K_S$.

Applying Lemma \ref{beta} once more, we may take a $\beta$
(arbitrarily near $\frac{4}{5}$) such that $\pi^*(K_X)|_S-\beta
\sigma^*(K_{S_0})$ is pseudo effective. Then $\xi\ge \frac{4}{5}$.
We have $p\ge 4$. Take $m=4$. Then
$\alpha=(4-1-\frac{1}{p}-\frac{1}{\beta})\xi\ge\frac{6}{5}>1$.
Theorem \ref{Key} gives $\xi\ge 1$.

So the only possibility is $\pi^*(K_X)|_S\cdot C=\xi=1$. We take the
family $\mathscr{C}$ on $X$ to be the set of images (via $\pi$)
of those fibers of $\Psi$. For a general member $C_0\in
\mathscr{C}$, one has $K_X\cdot C_0=\pi^*(K_X)\cdot C=1$. So the
case in Theorem \ref{partI} (ii) follows. This completes the proof
of Theorem \ref{partI}.
\end{setup}

\section{\bf Proof of the main theorem, Part II}

\begin{setup}{\bf Assumption.}\label{suppose}
Assume that a minimal threefold $X$ of general type is birationally
fibred by a family $\mathscr{C}$ of curves of geometric genus 2 with
$K_X\cdot C_0=1$ for a general member $C_0\in \mathscr{C}$. Suppose
$p_g(X)\ge 5$. We will show that $\varphi_4$ is not birational and
that the above family of curves is uniquely determined. This should
conclude Theorem \ref{main} by virtue of Theorem \ref{partI}.

After successive blow-ups, we may assume that $X'$ is the
birational model of $X$ which is fibred by genus two curves over
a base $W$. Denote by $f_{\mathscr{C}}:X'\longrightarrow W$ the
fibration. By assumption, a general fiber $C'$ of
$f_{\mathscr{C}}$ is a smooth curve of genus 2. Take a very
ample divisor $A$ on $W$. Set $H:=f^*_{\mathscr{C}}(A)$. Then
$|H|$ is a base point free linear system on $X'$.
\end{setup}

We now argue according to the value of $d = d(X) = \dim\varphi_1(X)$.

\begin{setup}\label{d=3} {\bf Claim. $d<3$.}\medskip

Suppose the contrary that $d=3$. Then $\dim\varphi_1(H)=2$ for a  general $H$. Note
that $|M_1|_H|\supset |M_1||_H$. We already know that
$\Phi_{|M_1||_H}$ is generically finite. So is $\Phi_{|M_1|_H|}$.
Thus, for a generic fiber $C'$ of $f_{\mathscr{C}}$ with $C'\subset
H$, $\dim\Phi_{|M_1|_H|}(C')=1$. The Riemann-Roch and Clifford's
theorem on $C'$ imply then that $(M_1\cdot C')_{X'}=(M_1|_H\cdot
C')_H\ge 2$. So $K_X\cdot C_0=\pi^*(K_X)\cdot C'\ge M_1\cdot C'\ge
2$, a contradiction.
\end{setup}

\begin{setup}\label{d=2}{\bf Claim. $\varphi_4$ is not birational whenever $d=2$.}

We consider a generic fiber $C'\subset H$ for a general $H$. The map
$\varphi_1|_H$ is exactly defined by the linear system
$|M_1||_H\subset |M_1|_H|$. If $\dim \varphi_1|_H(C')=1$, then the
Riemann-Roch and Clifford's theorem on $C'$ implies
$$K_X\cdot C_0=\pi^*(K_X)\cdot C'\ge M_1\cdot
C'=(M_1|_H\cdot C')_H\ge 2,$$ a contradiction. Thus, $\varphi_1$
maps a general $C'$ to a point. So Lemma 14 of \cite{KA} implies
that $f:X'\longrightarrow B$ birationally factors through
$f_{\mathscr{C}}$. Take further birational modifications to $X'$
such that $f$ factors through
$f_{\mathscr{C}} : X' \rightarrow W$ and a surjective morphism
$W \rightarrow B$ between normal projective surfaces. Then the uniqueness
of the Stein factorization says that $W \rightarrow B$ is birational and
we may assume $f=f_{\mathscr{C}}$ (after birationally modifying the base of $f_{\mathscr{C}}$).
This
simply means that the family $\mathscr{C}$ is exactly the
canonically induced family. Since $\pi^*(K_X)|_S\cdot C=K_X\cdot
C_0=1$, the argument in \ref{g2xi1} implies that $\varphi_4$ is
generically finite of degree 2.
\end{setup}

\begin{setup}\label{d=1}{\bf Claim. $\varphi_4$ is not birational
whenever $d=1$.}
\medskip

As we have shown above that $\dim\varphi_1|_H(C)=0$, $\varphi_1$
contracts a general fiber $C$ of $f_{\mathscr{C}}$ to a point.
Similarly we may suppose that $f: X'\longrightarrow B$ factors
through $f_{\mathscr{C}}:X'\longrightarrow W$ and
$\theta:W\longrightarrow B$.
\medskip

\begin{picture}(50,80) \put(100,0){$X$} \put(100,60){$X'$}
\put(170,0){$B$} \put(170,60){$W$}
\put(112,65){\vector(1,0){53}} \put(106,55){\vector(0,-1){41}}
\put(175,55){\vector(0,-1){43}} \put(114,58){\vector(1,-1){49}}
\multiput(112,2.6)(5,0){11}{-} \put(162,5){\vector(1,0){4}}
\put(133,70){$f_{\mathscr{C}}$} \put(180,30){$\theta$}
\put(92,30){$\pi$} \put(136,40){$f$}
\end{picture}
\bigskip

For a general point $b\in B$, the fiber $S = f^{-1}(b)$ has a
natural fibration $f_b:S\longrightarrow W_b = \theta^{-1}(b)$. A
general fiber of $f_b$ lies in the same numerical class as that
of $C$. Pick up a smooth fiber $C_b$ of $f_b$. By the
assumption, $\pi^*(K_X)|_S\cdot C_b=\pi^*(K_X)|_S\cdot
C=K_X\cdot C_0=1$. Since $p_g(X)\ge 5$, Lemma \ref{beta} tells
us that one may take a $\beta$ (arbitrarily near $\frac{4}{5}$)
such that $\pi^*(K_X)|_S-\beta \sigma^*(K_{S_0})$ is pseudo
effective, where $\sigma : S \rightarrow S_0$ is the smooth blow
down to a minimal model. So $\sigma^*(K_{S_0})\cdot C_b\le
\frac{5}{4}\pi^*(K_X)|_S\cdot C=\frac{5}{4}.$ This means that
$K_{S_0}\cdot L=\sigma^*(K_{S_0})\cdot C_b=1$ where
$L:=\sigma_*(C_b)$. This implies that $L^2>0$ and in fact
$K_{S_0}^2=L^2=1$ and $L \equiv K_{S_0}$ by the Hodge index
theorem. Now the surface theory tells us that $q(S_0)=0$ and
$1\le p_g(S_0)\le 2$ (by the Noether inequality). Note that
$g(W_b) \le q(S) = 0$, so $W_b$ is a smooth rational curve.
Hence $h^0(S, C_b)=2$. The computation on $S_0$ shows that
$p_a(L) = 2 = g(C_b)$. Thus $L \cong C_b$ is a smooth curve of
genus 2. If $S_0$ has the invariants $(K_{S_0}^2,
p_g(S_0))=(1,1)$ (resp. $(1,2)$), then the Neron Severi group of
$S_0$ has no torsion by Bombieri \cite{Bom} (resp. $S_0$ is
simply connected by \cite{BPV}). Therefore, $K_{S_0}\sim L$,
where the latter is movable. Hence $S_0$ must be of type (1,2).
This already shows that $\varphi_4$ is not birational. Since
$h^0(S_0, L) = p_g(S) = 2 = h^0(S, C_b)$, our $|C_b|$ is clearly
the movable part of $|\sigma^*(K_{S_0})|$. Since $C_b \equiv C$
and since $S$ is simply connected, we have $C_b \sim C$. Hence
the family of fibers $C$ of $f_{\mathscr{C}}$, when restricted
on $S$, is exactly the canonically induced family on $S$.
\end{setup}

\begin{setup}{\bf Conclusion to the main theorem.} Theorem \ref{partI},
\ref{d=3}, \ref{d=2} and \ref{d=1} imply Theorem \ref{main}. The
proof also shows that the family $\mathscr{C}$ in Theorem \ref{main}
is uniquely determined by $X$. We are done.
\end{setup}

\begin{setup}{\bf Proof of Theorem \ref{(1,2)}}
\begin{proof} We now prove Theorem \ref{(1,2)}.
Set $\tilde{b}:=g(\tilde{B})$.

The existence of $\tilde{f}$ implies that $\varphi_4$ is not
birational. Therefore $d=\dim\varphi_1(X)<3$ by \ref{d3}. The proof
of Theorem \ref{main} also implies the existence of the family
$\mathscr{C}$ of curves of genus 2. By taking further birational
modification, if necessary, we may suppose $X'=Z$. We also have the
birational modification $\pi=\nu:X'=Z\longrightarrow X$.
Denote by $\tilde{S}$ a general fiber of $\tilde{f}$.
Let $\tilde{\sigma}:\tilde{S}\longrightarrow \tilde{S_0}$ be the
blow-down onto the smooth minimal model.

{\bf Case 1}. $\tilde{b}>0$.

By the proof
of Lemma \ref{b>0}, we see that $\tilde{f}$ factors through $X$ and
thus $\nu^*(K_X)|_{\tilde{S}}\sim \tilde{\sigma}^*(K_{\tilde{S_0}})$.

We can actually get the family $\mathscr{C}$ by considering the
morphism $\tilde{\Psi}:Z\longrightarrow
\mathbb{P}(\tilde{f}_*\omega_{Z/\tilde{B}}^\checkmark).$ Let
$\tilde{C}\subset \tilde{S}$ be a general fiber of the induced
fibration by taking the Stein factorization of $\tilde{\Psi}$,
where $\tilde{S}$ is a smooth fiber of $\tilde{f}$. We have
$$\pi^*(K_X)\cdot \tilde{C}=\nu^*(K_X)\cdot
\tilde{C}=\tilde{\sigma}^*(K_{\tilde{S_0}})\cdot \tilde{C}=1.$$ So
the uniqueness of the family $\mathscr{C}$ says that $\mathscr{C}$
can be obtained by taking those images (via $\nu=\pi$) of all
fibers of $\tilde{\Psi}$.

{\bf Case 2}. $\tilde{b}=0$.

Let $\tilde{C}\subset \tilde{S}$ be a general fiber of the induced
fibration by taking the Stein factorization of $\tilde{\Psi}$,
where $\tilde{S}$ is a smooth fiber of $\tilde{f}$. Then $\tilde{C}$
is a smooth curve of genus 2. We want to show $\pi^*(K_X)\cdot
\tilde{C}=1$.

Considering the natural map
$$H^0(X', M_1)\longrightarrow H^0(\tilde{S}, M_1|_{\tilde{S}})\subset H^0(\tilde{S},
K_{\tilde{S}}),$$ we have $\pi^*(K_X)\ge M_1\ge 2\tilde{S}$ by the
fact that the vector space on the right has the dimension at most 2
and the assumption $p_g(X)\ge 5$. We take $J:=2\tilde{S}$ and will
apply Theorem \ref{Key}. Note that we have actually $\pi_J=\pi$ and
$f_J=\tilde{f}$. Clearly $p=2$. So Lemma \ref{beta} implies that one
can take a $\beta\mapsto \frac{2}{3}$ such that
$\pi^*(K_X)|_{\tilde{S}}-\beta\tilde{\sigma}^*(K_{\tilde{S_0}})$ is
pseudo-effective. Thus one has
$$\pi^*(K_X)\cdot \tilde{C}=\pi^*(K_X)|_{\tilde{S}}\cdot
\tilde{C}\ge \beta \tilde{\sigma}^*(K_{\tilde{S_0}})\cdot
\tilde{C}=\beta.$$ So $\tilde{\xi}=\pi^*(K_X)\cdot \tilde{C}\ge
\frac{2}{3}$. Now we can apply Theorem \ref{Key} to estimate
$\tilde{\xi}$. Take $m_1=5$. Then
$\alpha_1=(m_1-1-\frac{1}{p}-\frac{1}{\beta})\tilde{\xi}\ge
\frac{4}{3}.$ So Theorem \ref{Key} gives $\tilde{\xi}\ge
\frac{4}{5}$. Take $m_2=6$. Theorem \ref{Key} implies
$\tilde{\xi}\ge \frac{5}{6}$. In fact, an induction shows
$\tilde{\xi}\ge \frac{n-1}{n}$ for all $n\ge 8$. Thus
$\tilde{\xi}\ge 1$. Since $\tilde{\xi}\le
\tilde{\sigma}^*(K_{\tilde{S_0}})\cdot \tilde{C}=1$, we have
$\pi^*(K_X)\cdot \tilde{C}=1$.

The uniqueness of the family $\mathscr{C}$ means that
$\mathscr{C}$ can be obtained by taking those images (via
$\nu=\pi$) of all fibers of $\tilde{\Psi}$. We are done.
\end{proof}

\end{setup}

\section{\bf Application, new examples and open problems}

\begin{setup}{\bf Proof of Theorem \ref{app}}
\begin{proof}
(i) For a minimal threefold $X$ of general type with
$K_X^3<\frac{4}{3}(p_g(X)-2)$, one has $p_g(X)>2$ since $K_X^3>0$.
By assumption, we have $p_g(X)\ge 12$. We still consider the
canonical map and keep the same notation as above.

If $d=3$, then Kobayashi \cite{Kob} proved $K_X^3\ge 2p_g(X)-6\ge
\frac{4}{3}(p_g(X)-2)$, a contradiction.

If $d=2$ and $g(C)\ge 3$, we have shown in \ref{d2} that $\xi\ge
\frac{7}{4}$. So one has $K_X^3\ge
\frac{7}{4}(p_g(X)-2)>\frac{4}{3}(p_g(X)-2)$, also a contradiction.

If $d=2$, $g(C)=2$ and $\xi>1$, we have shown in \ref{g2xi>1} that
$\xi\ge \frac{5}{4}$ since $p_g(X)>5$. We may use Theorem \ref{Key}
to go on estimating $\xi$. Recall that we have $p=1$. Since
$p_g(X)\ge 12$, we may take $\beta\ge 10$. Take $m=3$. Then
$\alpha=(3-1-1-\frac{1}{\beta})\xi\ge \frac{9}{8}>1$. Theorem
\ref{Key} gives $\xi\ge \frac{4}{3}$. So $K_X^3\ge
\frac{4}{3}(p_g(X)-2)$, a contradiction.

If $d=1$ and $(K_{S_0}^2, p_g(S_0))\ne (1,2)$, then the results in
\cite{MRL} (Theorem 3.2(2), 3.4, 3.5) show that $K_X^3\ge
\frac{3}{2}p_g(X)-\frac{5}{2}>\frac{4}{3}(p_g(X)-2)$, a
contradiction.

Therefore, it is true that either $d=2$, $g(C)=2$ and $\xi=1$ (since
$p_g(X)>5$) or $d=1$ and $(K_{S_0}^2, p_g(S_0))=(1,2)$. Theorem
\ref{Key} implies that $\varphi_4$ is not birational.

(ii) By Theorem 4.1 in \cite{JMSJ3}, it is true that either $d=2$
and $g(C)=2$, or $d=1$ and $(K_{S_0}^2, p_g(S_0))=(1,2)$ whenever
$K_X^3<2p_g(X)-6$. We have to exclude the case with $d=2$ and
$\xi>1$. Noting that $\xi$ must be an integer, we have $\xi\ge 2$.
So, in that case, we still have $K_X^3\ge 2(p_g(X)-2)$, a
contradiction. Therefore, Theorem \ref{Key} implies that $\varphi_4$
is not birational.
\end{proof}
\end{setup}

\begin{exmp}
M. Kobayashi (see Proposition 3.2 in \cite{Kob}) has constructed a
family of canonically polarized smooth threefolds $Y$ satisfying the equality
$$K_Y^3=\frac{4}{3}p_g(Y)-\frac{10}{3}$$
where $p_g(Y)=7, 10, 13, \cdots$.

Theorem \ref{app} says that all examples above have non-birational
4-canonical maps.
\end{exmp}

The example below shows that the assumption $p_g(X) \ge 5$ in
Theorem \ref{main} is optimal.

\begin{exmp}\label{4} On ${\mathbb P}_{\mathbb C}^3$, take a smooth
hypersurface $S$ of degree 10. $S\sim 10H$ where $H$ is a
hyperplane. Let $\tau : X = Spec \oplus_{i=0}^1 \mathcal {O}(-5iH)
\rightarrow {\mathbb P}_{\mathbb C}^3$ be the double cover branched
along $S$. Then $X$ is a nonsingular canonical threefold with $K_X =
\tau^*H$, $K_X^3=2$ and $p_g(X)=4$ and $\Phi_{1}$ is a finite
morphism onto ${\mathbb P}^3$ of degree $2$. One may easily check
that $\Phi_{4}$ is also a finite morphism of degree $2$. Indeed, let
$C$ be the inverse on $X$ of a general line. Then $C$ is a
hyperelliptic curve of genus $4$ with $(4K_X)|_C = K_C + P + Q$,
where $\tau(P) = \tau(Q)$. Hence $\phi_4|_C$ is a degree 2 map (see
Section 6.5 of Iitaka's book \cite{Iitaka}).

It is not difficult to see that, for a generic irreducible curve
$C_0$ in any family of curves on $X$, $K_X\cdot C_0\ge 2$.
\end{exmp}

There are still several natural and unsolved problems:
\begin{setup}\label{open}{\bf Open problems.} Let $X$ be a minimal
projective threefold of general type with $\mathbb{Q}$FT singularities
and with $p_g(X)\ge 5$.
\begin{itemize}
\item[(1)]
Is it true that $\varphi_4$ is not birational if and only if $X$
is birationally fibred by a family of surfaces of type $(c_1^2, p_g)
= (1,2)$?
\item[(2)]
Is it possible to characterize the birationality of $\varphi_3$?
\item[(3)]
Is it true that $\dim \varphi_2(X)\ge 2$ when $p_g(X)$ is
bigger?
\end{itemize}
\end{setup}

We have only a partial answer to Problem \ref{open} (1). The above
problems might be very difficult, but very interesting. There has
not been any counter example to Open problems \ref{open}.
\bigskip

If $X$ is birationally fibred by surfaces of type $(1,2)$, then one
surely has the non-birationality of $\varphi_4$, which is hence a
natural condition in the result below.
\begin{thm}\label{ch(1,2)}
Let $X$ be a minimal projective threefold of general type with
${\bQ}FT$ singularities and with $p_g(X) \ge 5$. Suppose that
$\varphi_4$ is not birational. Then $X$ is birationally fibred by
surfaces of type $(1, 2)$ in the sense of Theorem $\ref{(1,2)}$ if and
only if either
\begin{itemize}
\item[(i)]
$\dim \Phi_{|K_X|}(X) = 1$, or
\item[(ii)]
$W' = \Phi_{|K_X|}(X)$ is of dimension $2$; further, with the
notation $g = s \circ f$ of Section $2$, $W' \subseteq
{\bP}^{p_g(X)-1}$ is covered by lines $\ell'$, these $\ell'$ move in
a pencil $\Lambda'$, and the pull back of the hyperplane $s^* H_{W'}
\ge \ell_1 + \ell_2$ for two generic irreducible members $\ell_i$ in
the movable part of $s^*\Lambda'$.
\end{itemize}
\end{thm}
\begin{proof}
By Theorem \ref{partI}, we may assume that $\dim \, W' = 2$. Suppose
that $X$ is fibred by surfaces of type $(1, 2)$ in the sense and
notation of Theorem \ref{(1,2)}. We may assume that $\nu : Z
\rightarrow X$ is identical to $\pi : X' \rightarrow X$ in notation
of Section 2, after further blow up. So there is a pencil $\tilde{f}
: X' \rightarrow \tilde{B}$ with a general fiber $F$ of type
$(1,2)$. We have $h^0(X', K_{X'} - F_1 - F_2) \ge h^0(X', K_{X'} -
F_1) - p_g(F) \ge p_g(X) - 2 p_g(F) \ge 1$, by considering the exact
sequences below for general fibers $F_1, F_2$ of $\tilde{f}$:
$$
0 \rightarrow \OO(K_{X'} - F_1) \rightarrow \OO(K_{X'}) \rightarrow
\OO(K_{F_1}) \rightarrow 0$$
$$0 \rightarrow \OO(K_{X'} - F_1 - F_2) \rightarrow \OO(K_{X'} - F_1)
\rightarrow \OO(K_{F_1}) \rightarrow 0.
$$
Thus $F_1 + F_2 \le M_1$ (the movable part of $K_{X'}$ or
$\pi^*K_X$). By Theorem \ref{partI}, the general fibers $C$ of $f :
X' \rightarrow B$ are curves of genus 2 with $\pi^*K_X \cdot C = 1$.
Since $C \cdot 2F \le C \cdot M_1 \le C \cdot \pi^*K_X = 1$, one has
$C \cdot F = 0$, so the general $C$ are contained in general fibers
$F$. By Lemma 14 in Kawamata \cite{KA} (after further blowing up
$X'$ and $B$ if necessary), $\tilde{f}$ factors as $\tilde{f} = \tau
\circ f$ where $\tau: B \rightarrow \tilde{B}$ is a surjective
morphism (from a surface to a curve) with connected fiber. Let $M_1
= f^*H_B$, $F_i = \tilde{f}^{-1}(\tilde{b}_i)$ and $\ell_i =
\tau^{-1}(\tilde{b}_i)$. Then $F_i = f^* \ell_i$ and $H_B \ge \ell_1
+ \ell_2$. Also $\ell_i = f(F_i)$ is a smooth rational curve because
$g(\ell_i) \le q(F_i)$ and $F_i$ is of type $(1, 2)$ (so $q(F_i) =
0$ by \cite{BPV}). Set $H_B = s^*H_{W'}$ with $H_{W'} = H' |_{W'}$
the restriction of a hyperplane $H'$ on ${\bP}^{p_g(X) - 1}$.

We assert that $H_B \cdot \ell_1 = 1$, or equivalently $M_1 \cdot
F_1 = C_1$ where $F_1$ is a general fiber of $f$. Indeed, $K_{F_1} =
K_{X'} |_{F_1} = (M_1 + E_1' + E_{\pi}) |_{F_1} \ge M_1 |_{F_1} =
f^*(H_B \cdot \ell_1) = \sum_{i=1}^e C_i$. Here $C_i$ are fibers of
the rational free pencil $F_1 \rightarrow \ell_1$, and $e \ge 1$
because $H_B$ is nef and big. Since $p_g(F_1) = 2$, we have $e = 1$.
The assertion is proved.

By the projection formula, $1 = s^* H_{W'} \cdot \ell_j = H_{W'} \cdot \ell_j' =
H' \cdot \ell_j'$,
where $\ell_j' := s_*(\ell_j)$. So $\ell_j'$ is irreducible and smooth; indeed it is
a line in ${\bP}^{p_g(X)-1}$; also
$\ell_j \rightarrow \ell_j'$ is birational and finite and hence an isomorphism.

Conversely, suppose that $H_B = s^*H_{W'} \ge \ell_1 + \ell_2$
where $\ell_j$ are generic irreducible members in a free pencil $\Lambda$,
the movable part of $s^* \Lambda'$ (on
the surface $B$; here $B$ and $X'$ are further blown up if necessary)
parametrized by a smooth curve $\tilde{B}$.
Then $M_1 = f^*H_B \ge F_1 + F_2$, where
the surface $F_j := f^*\ell_j$ is again parametrized by $\tilde{B}$.
We have only to show that $F_j$ is of type $(1, 2)$.
Let $\sigma : F = F_j \rightarrow F_0$ be the smooth blow down to a minimal model.

We assert that $C \cdot \sigma^* F_0 = 1$ with $C$
a general fiber of $f|_{F_j} : F_j \rightarrow \ell_j = f(F_j)$
(and also a general fiber of $f$).
By Theorem \ref{partI}, $C \cdot \pi^*K_X = 1$.
If the parametrizing curve $\tilde{B}$ of surfaces $F$
has genus $\ge 1$, then $\sigma^* K_{F_0} = \pi^*K_X |_{F}$ by the proof of
Lemma \ref{b>0}. So the assertion is clear.
If $g(\tilde{B}) = 0$, then by Lemma \ref{beta}, we have
$1 \le \sigma^*K_{F_0} \cdot C \le \frac{3}{2} \pi^* K_X \cdot C = 3/2$.
Hence the assertion is true. Now the assertion and the argument in
Claim \ref{d=1} imply that $F$ is of type (1, 2). This proves the theorem.
\end{proof}

\noindent{\bf Acknowledgment.} The first author would like to thank
the National University of Singapore for the support during his
visit in the Spring of 2006. Both authors would like to thank
Fabrizio Catanese for his query on the possibility of characterizing
the birationality of $\varphi_4$ of threefolds of general type.
Thanks are also due to Eckart Viehweg who spent a lot of time to
read the paper and to give us valuable comments.
The authors like to thank the referee for very careful reading
and valuable and detailed suggestions.


\end{document}